\documentclass[11pt]{article}

\usepackage{amsmath,amssymb,amsthm}
\usepackage{graphicx}
\usepackage{geometry}
\usepackage{setspace}
\usepackage{hyperref}
\usepackage{float}
\usepackage{bm}

\newtheorem{remark}{Remark}

\title{Delay-Induced Stability Transitions in Directed Signed Consensus Networks}

\author{
Hui Wu\\
Department of Mathematics\\
Clark Atlanta University\\
\texttt{hwu@cau.edu}
}

\date{\today}

\begin{document}

\maketitle

\begin{abstract}
We study delay-induced transitions in consensus dynamics on signed networks with a ring topology. The proposed model is formulated as a system of delay differential equations incorporating both cooperative and antagonistic interactions, as well as heterogeneous time delays.

We perform a stability analysis by deriving the associated characteristic equation and examining the real parts of its eigenvalues. It is shown that the stability of the consensus state depends critically on the magnitude of the delays. In particular, increasing time delays may destabilize the system and induce transitions from consensus to bounded non-convergent behavior or even instability.

A phase diagram in the parameter space is constructed to identify different dynamical regimes. Numerical simulations are presented to validate the theoretical results and to illustrate the delay-induced transitions.

Such delay-induced transitions have also been reported in various biological and engineered systems, highlighting the universal role of time delays in shaping collective dynamics.
\end{abstract}

\noindent\textbf{Keywords:} time delay; delay differential equations; consensus dynamics; signed networks; dynamical systems; stability analysis

\section{Introduction}

Time-delay effects play a fundamental role in the dynamics of complex systems, where finite signal transmission speeds, processing times, and feedback mechanisms inevitably introduce delays into interactions among system components \cite{hale,michiels,erneux,stepan}. Such delays are well known to significantly affect system stability and may induce oscillations, bifurcations, or even instability in otherwise well-behaved systems.

Consensus problems in networked multi-agent systems have attracted considerable attention due to their wide range of applications in engineering, biology, and social systems \cite{olfati2007,ren2005,jadbabaie2003,moreau2005}. In the classical consensus framework, agents update their states based on local interactions, typically modeled through diffusive coupling. It has been shown that time delays can degrade performance or even destabilize consensus dynamics \cite{olfati_delay,frasca2008,seuret2015}.

In many real-world systems, interactions are not purely cooperative. Instead, antagonistic or competitive effects may arise, leading to the concept of signed networks \cite{altafini2013,proskurnikov2017,shi2013}. In such networks, positive and negative couplings coexist and give rise to more complex dynamical behaviors, including clustering, oscillation, and instability.

The interplay between time delays and signed interactions introduces additional challenges in understanding system behavior. Delay-induced transitions have been observed in various biological and engineered systems, highlighting the universal role of time delays in shaping collective dynamics \cite{erneux,stepan,gu2003}. Similar phase-delay mechanisms appear in optical ring resonators, where propagation delays lead to phase accumulation and determine resonance conditions \cite{zhang1988}. These observations suggest that delays can fundamentally alter coordination patterns in networked systems.

In this paper, we study a delayed consensus model with cooperative and antagonistic couplings on a directed ring network. The choice of a ring topology is motivated by both analytical tractability and its role as a canonical network structure. Ring networks represent one of the simplest nontrivial directed graphs while capturing essential features of cyclic information flow and feedback. Moreover, their circulant structure allows the application of Fourier mode decomposition, which decouples the dynamics into independent modes and facilitates stability analysis \cite{olfati2007,dorfler2014}. For these reasons, ring networks are widely used as benchmark models in the study of network dynamics and synchronization.

The form of the coupling in our model follows the standard diffusive interaction framework commonly used in consensus and synchronization problems \cite{ren2005,jadbabaie2003}. In particular, terms of the form $(x_j(t-\tau)-x_i(t))$ represent the tendency of each agent to adjust its state based on discrepancies with its neighbors. The introduction of two distinct delays reflects the fact that information transmission in different directions or channels may occur at different rates, which is often observed in practical systems \cite{michiels,erneux}. The inclusion of antagonistic interactions further enriches the dynamics and may lead to nontrivial behaviors such as clustering, oscillation, or instability.

By combining diffusive coupling, signed interactions, and heterogeneous delays within a ring structure, the proposed model provides a minimal yet nontrivial framework for studying delay-induced stability transitions in networked systems. Our main objective is to characterize how time delays influence the stability of the consensus state and to identify the resulting dynamical regimes.

The main contributions of this paper are summarized as follows:
\begin{itemize}
    \item We propose a delayed consensus-type model on a directed signed ring network with distinct delays in cooperative and antagonistic couplings.
    \item By exploiting the circulant structure, we reduce the system to scalar delay differential equations associated with Fourier modes.
    \item We derive the corresponding characteristic equation and analyze the stability of the consensus state.
    \item Numerical simulations reveal a phase diagram in the $(\tau_1,\tau_2)$-plane with three distinct regimes: consensus, bounded non-convergent behavior, and instability.
\end{itemize}

\section{Model formulation}

Consensus dynamics in networked systems are commonly modeled through diffusive coupling between neighboring agents \cite{olfati2007,ren2005,jadbabaie2003}. In many real-world systems, however, interactions are not purely cooperative. Instead, antagonistic or competitive effects may arise, leading to the notion of signed networks \cite{altafini2013,proskurnikov2017}.

In addition, communication and processing delays are inevitable in both engineered and natural systems, and they can significantly influence system stability and collective behavior \cite{hale,michiels,erneux}. Motivated by these considerations, we consider a delayed consensus model on a directed signed network, where each node interacts with its neighbors through both cooperative and antagonistic couplings, subject to distinct time delays.

Specifically, we study a ring-type directed network with $N$ agents, where each agent receives delayed information from its forward and backward neighbors. The system is governed by the following delay differential equations:
\begin{equation}
\dot{x}_i(t) = K_p \big(x_{i+1}(t-\tau_1) - x_i(t)\big) 
- K_n \big(x_{i-1}(t-\tau_2) - x_i(t)\big),
\quad i = 1,\dots,N,
\label{eq:model_scalar}
\end{equation}
with periodic boundary conditions.

Here, $K_p > 0$ represents the cooperative coupling strength, while $K_n > 0$ models antagonistic interaction. The delays $\tau_1$ and $\tau_2$ correspond to forward and backward information transmission times, respectively.

Rewriting the system, we obtain
\begin{equation}
\dot{x}_i(t)
= -(K_p - K_n)x_i(t)
+ K_p x_{i+1}(t-\tau_1)
- K_n x_{i-1}(t-\tau_2),
\label{eq:model_expanded}
\end{equation}
which highlights the balance between self-dynamics and delayed interactions.

Such delay-coupled network models have been widely used in the study of multi-agent coordination, synchronization, and distributed control systems \cite{olfati_delay,frasca2008,seuret2015}. The introduction of signed interactions further enriches the dynamics and may lead to nontrivial behaviors such as clustering, oscillation, or instability, as will be illustrated in the subsequent sections.

\section{Stability analysis}

To analyze the stability of the consensus state, we exploit the circulant structure of the network and apply Fourier mode decomposition \cite{olfati2007,dorfler2014}. Each mode evolves independently according to
\begin{equation}
\dot{z}_k(t)
= -(K_p - K_n) z_k(t)
+ K_p e^{i\theta_k} z_k(t-\tau_1)
- K_n e^{-i\theta_k} z_k(t-\tau_2),
\label{eq:mode_eq}
\end{equation}
where $\theta_k = \frac{2\pi k}{N}$.

Assuming solutions of the form $z_k(t)=e^{\lambda t}$, we obtain the characteristic equation \cite{michiels,niculescu2001}
\begin{equation}
\lambda + (K_p - K_n)
- K_p e^{-\lambda \tau_1} e^{i\theta_k}
+ K_n e^{-\lambda \tau_2} e^{-i\theta_k}
= 0.
\label{eq:char}
\end{equation}

The roots $\lambda \in \mathbb{C}$ of \eqref{eq:char} determine the stability of each mode. In general, these roots may be complex due to the presence of time delays.

\medskip
\noindent\textbf{Stability condition.}
The consensus state is asymptotically stable if and only if all roots of \eqref{eq:char} satisfy
\begin{equation}
\mathrm{Re}(\lambda)<0
\qquad \text{for all modes } k\neq 0.
\end{equation}
The mode $k=0$ corresponds to perturbations along the consensus manifold and therefore does not affect transverse stability.

\medskip
\noindent\textbf{Stability boundary.}
To determine the stability boundary, we consider purely imaginary roots $\lambda=i\omega$, $\omega\in\mathbb{R}$ \cite{erneux}. Substituting into \eqref{eq:char} yields
\begin{equation}
i\omega + (K_p - K_n)
= K_p e^{-i\omega \tau_1} e^{i\theta_k}
- K_n e^{-i\omega \tau_2} e^{-i\theta_k}.
\end{equation}

Separating real and imaginary parts, we obtain
\begin{align}
K_p - K_n
&= K_p \cos(\theta_k - \omega \tau_1)
- K_n \cos(\theta_k + \omega \tau_2),
\label{eq:realpart}\\
\omega
&= K_p \sin(\theta_k - \omega \tau_1)
+ K_n \sin(\theta_k + \omega \tau_2).
\label{eq:imagpart}
\end{align}
These equations implicitly define the stability switching curves in the $(\tau_1,\tau_2)$-parameter space.

\begin{remark}
Due to the transcendental nature of the characteristic equation, explicit expressions for the roots are generally unavailable. Therefore, numerical simulations are employed to explore the stability regions and validate the analytical predictions.
\end{remark}

\section{Numerical results}

In this section, we present numerical simulations to illustrate the theoretical analysis and to explore the influence of time delays on the collective dynamics of the system.

All simulations are performed with a network size $N=20$. We have verified that the qualitative behavior of the system remains consistent for different network sizes, indicating that the observed dynamical regimes are not sensitive to the choice of $N$.

\subsection{Phase diagram in the $(\tau_1,\tau_2)$-plane}

We first investigate the global behavior of the system in the delay parameter space. For a fixed network size and coupling strengths, the system is simulated over a grid of delay values $(\tau_1,\tau_2)$. For each parameter pair, the long-term behavior is classified based on the consensus error
\begin{equation}
V(t) = \frac{1}{N}\sum_{i=1}^N \left(x_i(t)-\bar{x}(t)\right)^2,
\qquad
\bar{x}(t)=\frac{1}{N}\sum_{i=1}^N x_i(t).
\label{eq:Vt}
\end{equation}

\begin{figure}[H]
\centering
\includegraphics[width=0.72\textwidth]{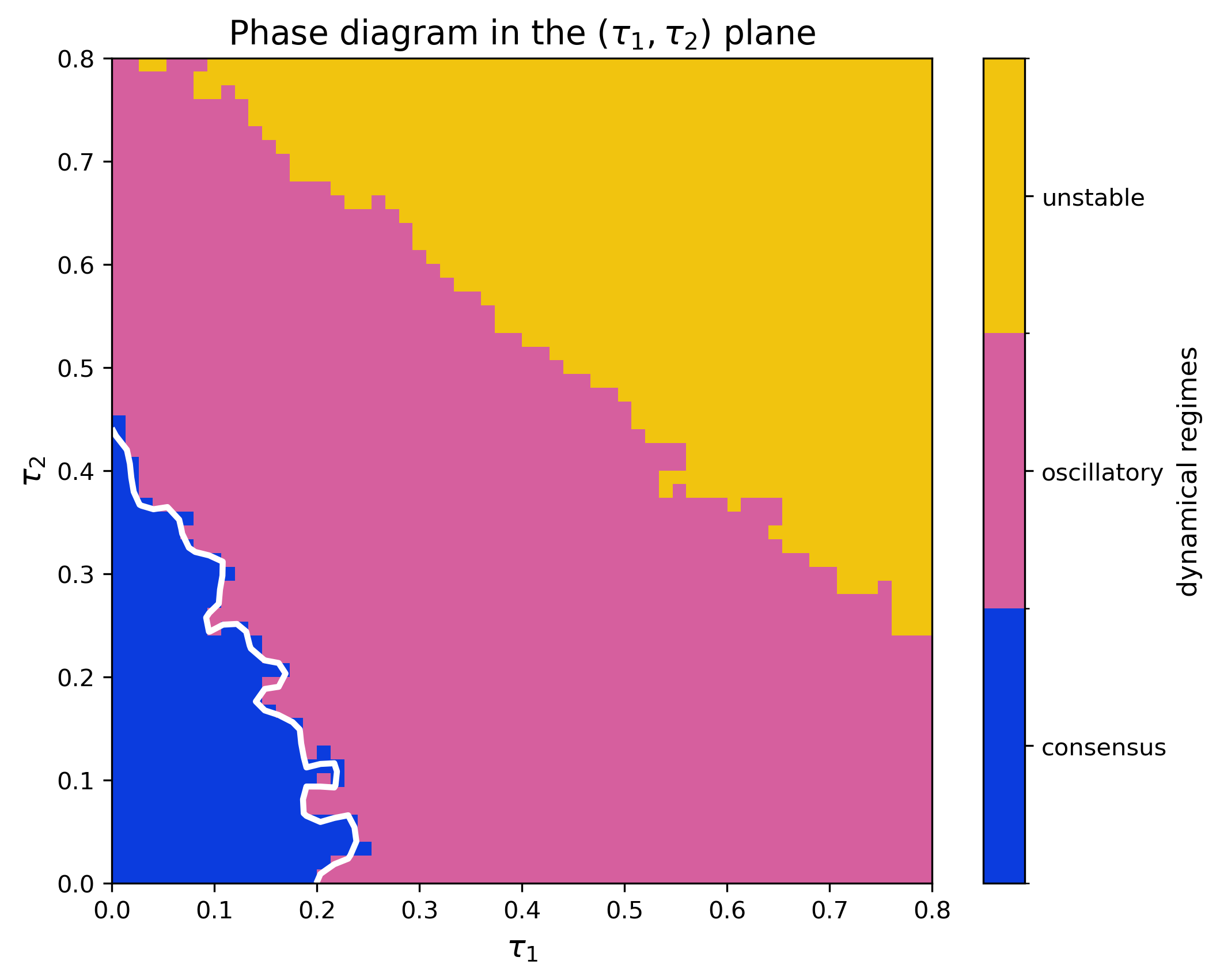}
\caption{Phase diagram of the delayed system in the $(\tau_1,\tau_2)$ plane. Three dynamical regimes are observed: consensus, bounded non-convergent behavior, and instability. The white curve indicates the numerically detected transition boundary between the consensus and bounded regimes.}
\label{fig:phase}
\end{figure}

Figure~\ref{fig:phase} shows the resulting phase diagram in the $(\tau_1,\tau_2)$-plane. Three distinct dynamical regimes can be clearly identified:
\begin{itemize}
    \item \textbf{Consensus region}: the consensus error $V(t)$ decays to zero, indicating that all nodes converge to a common state;
    \item \textbf{Bounded region}: the consensus error remains bounded but does not converge to zero;
    \item \textbf{Instability region}: the consensus error grows rapidly, indicating loss of coordinated behavior.
\end{itemize}

The white curve in Figure~\ref{fig:phase} represents the numerically detected transition boundary between the consensus and bounded regimes. It can be observed that the stability region shrinks as the delays increase, indicating that time delays have a destabilizing effect on the system.

\subsection{Time evolution of the consensus error}

To further illustrate the dynamical behavior in each regime, we select three representative parameter pairs from the phase diagram and plot the time evolution of the consensus error $V(t)$.

\begin{figure}[H]
\centering
\includegraphics[width=0.72\textwidth]{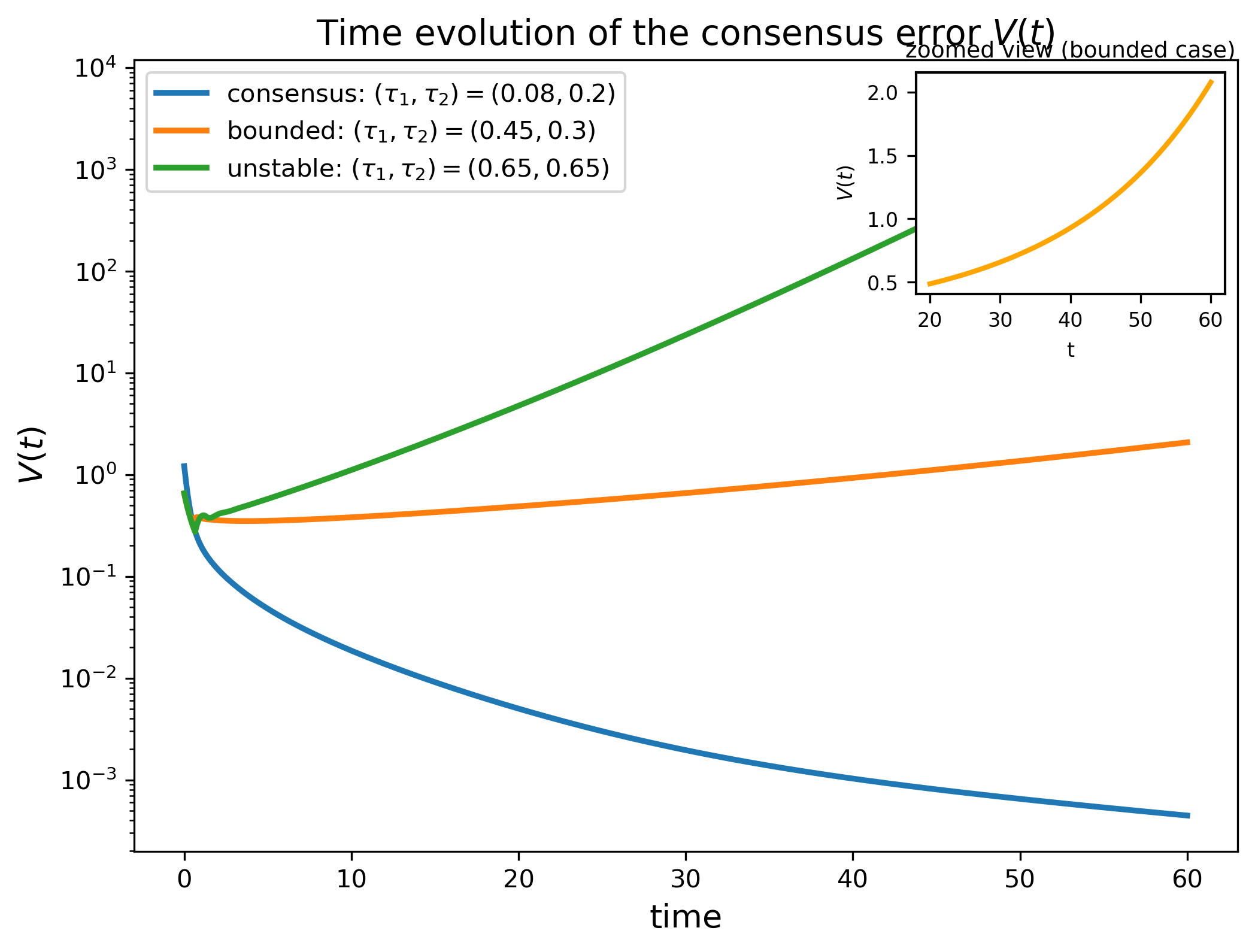}
\caption{Time evolution of the consensus error $V(t)$ for three representative parameter pairs. A logarithmic scale is used in the main panel. For small delays, $V(t)$ decays to zero, indicating consensus. For intermediate delays, $V(t)$ remains bounded but does not converge. For larger delays, $V(t)$ grows exponentially, indicating instability. The inset shows a zoomed view of the bounded non-convergent case.}
\label{fig:vt}
\end{figure}

As shown in Figure~\ref{fig:vt}, for small delays, the consensus error decays monotonically to zero, confirming asymptotic consensus. For intermediate delays, the consensus error remains bounded but does not converge, indicating the loss of global agreement. For larger delays, the consensus error grows rapidly, suggesting that the system becomes unstable.

A logarithmic scale is used in the main panel of Figure~\ref{fig:vt} to clearly distinguish the different dynamical behaviors. The inset shows a zoomed view of the bounded case, highlighting that the consensus error does not converge even over long time intervals.

In addition to the long-term behavior, we observe a short transient phase at early times, during which the trajectories rapidly adjust from the initial conditions. This transient behavior reflects the interaction between instantaneous self-dynamics and delayed coupling terms before the delayed feedback becomes dominant.

The logarithmic scale in Figure~\ref{fig:vt} allows a clear comparison of growth rates across regimes. In the consensus case, the consensus error decays approximately exponentially, indicating that all transverse modes are stable. In the bounded case, the error remains finite and exhibits slow variation over time, suggesting that the dominant modes have real parts close to zero. In contrast, the instability case shows a clear exponential growth, corresponding to eigenvalues with positive real parts.

The inset further highlights the bounded case by zooming into the interval $t\in[20,60]$, confirming that the consensus error does not converge to zero even over long time horizons. This excludes the possibility of slow convergence and supports the classification as a non-convergent bounded regime.

These numerical observations are consistent with the stability analysis, where the system behavior is determined by the location of the roots of the characteristic equation. In particular, transitions between regimes correspond to changes in the sign of the real part of the dominant eigenvalues.

\subsection{Discussion}

The numerical results reveal a clear transition in the system behavior as the delays increase. In particular, increasing delays can drive the system from a stable consensus state to bounded non-convergent dynamics and eventually to instability. This transition is consistent with the stability switching mechanism predicted by the characteristic equation.

These observations suggest that time delays play a crucial role in shaping the collective dynamics of the network. In applications such as neural systems, where transmission delays are unavoidable and may increase with aging, such delay-induced transitions may provide a possible mechanism for the degradation of coordinated activity.

\section{Conclusion}

In this work, we investigated the impact of time delays on consensus dynamics in directed signed networks with cooperative--antagonistic interactions. By combining analytical derivations and numerical simulations, we identified three distinct dynamical regimes: consensus, bounded non-convergent behavior, and instability.

Our theoretical analysis, based on modal decomposition and characteristic equations, provides insight into how delays affect the stability of each mode. The numerical phase diagram further illustrates the existence of stability regions and transition boundaries in the $(\tau_1,\tau_2)$-parameter space.

In particular, the results show that increasing delays can induce qualitative transitions in the system dynamics, leading to loss of consensus and eventual instability. The time evolution of the consensus error confirms these behaviors and highlights the role of delays in shaping long-term dynamics.

Such delay-induced transitions have also been reported in various biological and engineered systems, highlighting the universal role of time delays in shaping collective dynamics.

Future work may focus on extending the analysis to heterogeneous delays, larger-scale networks, and more general interaction topologies. In addition, developing sharper analytical conditions for stability boundaries remains an interesting direction for further research.

\end{document}